\documentstyle[12pt]{article}
\def\ep{\epsilon}  \def\concat{{\;\hat{}\;}}
\def\om{\omega}    \def\ka{\kappa}   \def\ga{\gamma}  \def\Ga{\Gamma}
\def\al{\alpha}    \def\be{\beta}    \def\de{\delta}  \def\rh{\rho}
\def\De{\Delta}    \def\Si{\Sigma}   \def\forces{\mid\vdash}
\def\su{\subseteq} \def\ra{\rangle}  \def\la{\langle} \def\se{\setminus}
\def\supp{{\rm supp}}  \def\dom{{\rm dom}}    \def\cof{{\rm cof}}
\def\rmiff{{\mbox{ iff }}}  \def\rmand{{\mbox{ and }}}
\def\rmor{{\mbox{ or }}}
\def\pf{\par\noindent{\bf Proof:}\par}
\def\qed{\par\noindent$\Box$\par\medskip}
\def\pfof#1{{\bigskip \noindent{\bf  Proof of Theorem \ref{#1}.}\medskip  }}

\newtheorem{theorem}{Theorem} \newtheorem{lemma}[theorem]{Lemma}

\def\res{\upharpoonright}
\def\rr{{\Bbb R}} \def\rat{{\Bbb Q}}
\def\po{{\Bbb P}} \def\cc{{\goth c}}
\input amssym.def
\input amssym

\begin{document}

\begin{center}
  {\Large Orthogonal Families of Real Sequences}
\end{center}

\begin{center}
  Arnold W. Miller\\
   and \\
  Juris Steprans\\
\end{center}

\bigskip

For $x,y\in\rr^\om$ define the inner product
$$(x,y)=\Sigma_{n\in\om}x(n)y(n)$$
which may not be finite
or even exist. We say that $x$ and $y$ are orthogonal if
$(x,y)$ converges and equals $0$.

Define $l_p$ to be the set of all $x\in\rr^\om$ such that
$$\sum_{n\in\om}|x(n)|^p<\infty.$$

For Hilbert space, $l_2$, any family of pairwise
orthogonal sequences must be countable.
For a good introduction
to Hilbert space, see Retherford \cite{R}.

\begin{theorem} \label{thm1}
There exists a pairwise orthogonal family $F$ of size continuum
such that $F$ is a subset of $l_p$ for every $p>2$.
\end{theorem}

It was already known\footnote{Probably Kunen was the first.
His example is given in the proof of Theorem \ref{thm1-5}.
Earlier work was done by Abian and examples
were also constructed by Keisler and Zapletal independently.
At any rate, we know definitely that we didn't do it first.}
that there exists a family
of continuum many pairwise orthogonal elements of $\rr^\om$.
A family $F\subseteq\rr^\om\se {\bf 0}$ of pairwise
orthogonal sequences is
orthogonally complete or a maximal orthogonal family
iff the only element of $\rr^\om$ orthogonal to
every element of $F$ is ${\bf 0}$, the constant 0 sequence.

It is somewhat surprising that Kunen's
perfect set of orthogonal elements is maximal (a fact first
asserted by Abian).  MAD families, nonprincipal ultrafilters,
and many other such maximal objects cannot be even Borel.

\begin{theorem} \label{thm1-5}
There exists a perfect maximal orthogonal family of
elements of $\rr^\om$.
\end{theorem}

Abian raised the question of what are the
possible cardinalities of maximal orthogonal families.

\begin{theorem} \label{thm2}
In the Cohen real model there is a maximal orthogonal
set in $\rr^\om$ of cardinality $\om_1$, but there is no maximal
orthogonal set of cardinality $\ka$ with $\om_1<\ka<\cc$.
\end{theorem}

By the Cohen real model we mean any model obtained by forcing
with finite partial functions from $\ga$ to $2$, where the
ground model satisfies GCH and $\ga^\om=\ga$.

\begin{theorem} \label{thm3}
For any countable standard model $M$ of ZFC and
$\ka$ in $M$ such that $M\models \ka^\om=\ka$, there exists
a ccc generic extension $M[G]$ such that
the continuum of $M[G]$ is $\ka$ and in $M[G]$
for every infinite cardinal $\al\leq\ka$ there is a
maximal orthogonal family of cardinality $\al$.
\end{theorem}

\begin{theorem} \label{ma} (MA$_\ka$($\sigma$-centered))
Suppose $X\su \rr^\om$, $||X||\leq\ka$, $X\cap l_2$ is finite,
and for
every distinct pair $x,y\in X$ the inner product $(x,y)$ converges.  Then
there exists a $z\in \rr^\om\se l_2$ such that $z$ is orthogonal to every
element of $X$.
\end{theorem}

The question arises of whether uncountable families of
pairwise orthogonal elements of $\rr^\om$ must somehow determine
an almost disjoint family of subsets of $\om$.  In the
following result we give a perfect family of orthogonal
elements of $\rr^\om$ each of which has full support.
It is possible to modify Kunen's example, Theorem \ref{thm1-5},
or the method of Theorem \ref{thm1} to produce such
a perfect set by replacing the zeros by
a very small sequence of positive weights.
If this\footnote{
Briefly, use the  weights $a_n={1\over{2^{4n}}}$ off
of the comb on level $n$.
Use weights $+-b_n$ on the comb with
$b_0$ on the root node, $-b_n$ on the teeth, and $+b_n$
on the branch.
By choosing $b_{n+1}$ so that
$$b_0^2 + 2b_1^2 + ... + 2b_n^2 - 2b_{n+1}^2 +
\sum_{i\leq n+1} a_i^2(2^i-2) + \sum_{i>n+1} a_i^2 (2^i -4) =0,$$
the inner products will be zero.}
is done, then the resulting elements will be
``$l_2$ almost disjoint'' in the following sense.
Given $x$ and $y$ there will be $X$ and $Y$ almost
disjoint subsets of $\om$ such that
$$\sum_{n\notin X}x(n)^2<\infty \rmand \sum_{n\notin Y}y(n)^2<\infty.$$
Equivalently,
$$\sum_{n\in\om}\min\{|x(n)|,|y(n)|\}^2<\infty.$$
Note that the supports of $x$ and $y$ are
almost disjoint iff the function
$$n\mapsto\min\{|x(n)|,|y(n)|\}$$
is eventually zero.
Consequently, the
minimum function is one measure of the almost disjointedness of
$x$ and $y$.  Note, however that if the inner product
of $x$ and $y$ converges, then $\min\{|x(n)|,|y(n)|\}\to 0$
as $n\to\infty$.

\begin{theorem}\label{support}
There exists a perfect set $P\su\rr^\om$ such that
every pair of elements of $P$ are orthogonal,
$\supp(x)=\om$ for every $x\in P$, and if
we define
$$h(n)=\min\{|x(n)|:x\in P\}$$
then for every $p$
$$\sum_{n<\om}h(n)^p=\infty.$$
\end{theorem}

K.P.Hart raised the question of whether there could be a maximal
orthogonal family in $l_2$ which was not maximal in $\rr^\om$.
This was answered by Kunen and Steprans independently.

\begin{theorem} \label{l2}
(a) There exists $X$ which is a
maximal orthogonal family in $l_2$ such that
for all $n$ with
$1\leq n \leq\om$ there exists $Y\su \rr^\om\se l_2$ with $||Y||=n$
and
$X\cup Y$ a maximal orthogonal family in $\rr^\om$.  Furthermore,
every maximal orthogonal family containing $X$ is countable.

(b) There exists a perfect maximal orthogonal
family $P\su\rr^\om$ such that $P\cap l_2$ is a maximal
orthogonal family in $l_2$.
\end{theorem}

\begin{center}
  The Proofs
\end{center}

\pfof{thm1}

Here is the basic idea.
Take the full binary tree and attach pairwise disjoint finite sets $F_s$
to each node, see figure \ref{comb}.

Instead of taking branches, take "combs", which are
a branch together with nodes which are just off the branch,
e.g. for rightmost branch the comb would be as in figure \ref{comb}.
The elements of the comb will be the support of sequence.
Attach to the branch nodes $F_{111..1}$ positive weights and to the
off branch nodes $F_{111...10}$ negative weights.  Then when
two combs eventually disagree the lowest level pair of nodes gives
a negative value.  By choosing the
sizes of the $F_s$'s and attached weights correctly this negative
value will cancel out all the positive values above.

\begin{figure}
\unitlength=1.00mm
\special{em:linewidth 0.4pt}
\linethickness{0.4pt}
\begin{picture}(130.00,50.00)
\put(40.00,50.00){\makebox(0,0)[cc]{$F_{\la\ra}$}}
\put(20.00,30.00){\makebox(0,0)[cc]{$F_{0}$}}
\put(60.00,30.00){\makebox(0,0)[cc]{$F_{1}$}}
\put(10.00,10.00){\makebox(0,0)[cc]{$F_{00}$}}
\put(30.00,10.00){\makebox(0,0)[cc]{$F_{01}$}}
\put(50.00,10.00){\makebox(0,0)[cc]{$F_{10}$}}
\put(70.00,10.00){\makebox(0,0)[cc]{$F_{11}$}}
\put(25.00,35.00){\line(1,1){10.00}}
\put(55.00,35.00){\line(-1,1){10.00}}
\put(10.00,15.00){\line(1,1){10.00}}
\put(20.00,25.00){\line(1,-1){10.00}}
\put(50.00,15.00){\line(1,1){10.00}}
\put(60.00,25.00){\line(1,-1){10.00}}
\put(90.00,50.00){\makebox(0,0)[cc]{$F_{\la\ra}$}}
\put(80.00,40.00){\makebox(0,0)[cc]{$F_0$}}
\put(100.00,40.00){\makebox(0,0)[cc]{$F_1$}}
\put(90.00,30.00){\makebox(0,0)[cc]{$F_{10}$}}
\put(110.00,30.00){\makebox(0,0)[cc]{$F_{11}$}}
\put(100.00,20.00){\makebox(0,0)[cc]{$F_{110}$}}
\put(120.00,20.00){\makebox(0,0)[cc]{$F_{111}$}}
\put(110.00,10.00){\makebox(0,0)[cc]{$F_{1110}$}}
\put(130.00,10.00){\makebox(0,0)[cc]{$F_{1111}$}}
\put(83.00,43.00){\line(1,1){4.00}}
\put(97.00,43.00){\line(-1,1){4.00}}
\put(93.00,33.00){\line(1,1){4.00}}
\put(103.00,37.00){\line(1,-1){4.00}}
\put(103.00,23.00){\line(1,1){4.00}}
\put(113.00,27.00){\line(1,-1){4.00}}
\put(113.00,13.00){\line(1,1){4.00}}
\put(124.00,17.00){\line(1,-1){4.00}}
\end{picture}
\caption{Combs \label{comb}}
\end{figure}

Define the sequence $r_n$ by $r_0=1$ and
$$r_{n+1}=\sum_{i \leq n} r_i.$$
(This makes $r_n=2^{n-1}$ for $n>0$, but it is irrelevant.)
Let $p_n>2$ be a sequence decreasing to 2.
Next construct integers $k_n>r_n$ and reals $\ep_n>0$ such
that
$$ \ep_n^2 \cdot k_n =r_n$$
and
$$ \ep_n^{p_n}\cdot k_n \leq {1\over n^2}.$$
To do this first pick $k_n$ so that
$$ k_n^{{p_n\over 2}-1}\geq n^2\cdot r_n^{{p_n}\over 2}$$
then let
$$\ep_n^2={r_n\over k_n}.$$
Thus
$$\ep_n^{p_n}\cdot k_n =({r_n\over k_n})^{p_n\over 2}\cdot k_n=
{{r_n^{p_n/2}}\over{k_n^{p_n/2-1}}}\leq {1\over n^2}.$$

Let $2^\om$ be the set of infinite sequences of 0 and 1's and
let $2^{<\om}$ be the set of finite sequences of 0 and 1's.
Let $\{F_s: s\in 2^{<\om}\}$ be pairwise disjoint subsets of
$\om$ such that $||F_s||=k_n$ when $s\in 2^n$.

For $x\in 2^\om$ and $n\in\om$ define
$$s^x_{n+1}=x |_{n+1}$$
and
$$t^x_{n+1}=x |_n \concat (1-x(n)).$$
Thus $s^x_{n+1}$ is the first $n+1$ bits of $x$ while
$t^x_{n+1}$ is the first $n$ bits of $x$ followed by the opposite bit $1-x(n)$.
Define $y_x\in\rr^\om$ by
$$y_x(m)=\left\{\begin{array}{rl}
          \ep_{n+1} & \mbox{if $m\in F_{s^x_{n+1}}$ for some $n$}\\
         -\ep_{n+1} & \mbox{if $m\in F_{t^x_{n+1}}$ for some $n$}\\
\sqrt{2}\cdot\ep_0  & \mbox{ if $m\in F_{\la\ra}$ } \\
           0        & \mbox{ otherwise }
         \end{array}\right.$$

First note that $y_x\in l_p$ for any $p>2$:
$$ \sum |y_x(m)|^p=(\sqrt{2}\ep_0)^p\cdot k_0 +\sum (\ep_{n+1}^p
||F_{s_{n+1}}||+
     \ep_{n+1}^p ||F_{t_{n+1}}||)$$
and
$$ \sum (\ep_{n+1}^p ||F_{s_{n+1}}||+
                           \ep_{n+1}^p ||F_{t_{n+1}}||)=
                   \sum 2\ep_{n+1}^p k_{n+1}.$$

But for all but finitely many $n$ we have that $p_n<p$ and
so $\ep_n^p k_n < {1\over n^2}$.

Now we see that for distinct $x,x^\prime\in 2^\om$
that $y_x$ and $y_{x^\prime}$ are orthogonal.
Take $N<\om$ so that $x|_N=x^\prime|_N$ but $x(N)\not=x^\prime(N)$.
Thus $s_n^x=s_n^{x^\prime}$ and $t_n^x=t_n^{x^\prime}$ for $n\leq N$,
but $s_{N+1}^x=t_{N+1}^{x^\prime}$ and $t_{N+1}^x=s_{N+1}^{x^\prime}$.
Therefore
$$(y_x,y_{x^\prime})=(\sum_{n\leq N}2\ep_n^2k_n)- 2\ep_{N+1}^2k_{N+1} =
(\sum_{n\leq N}2r_n)- 2r_{N+1}=0.$$

\pfof{thm1-5}

\begin{figure}
\unitlength=.9mm
\special{em:linewidth 0.4pt}
\linethickness{0.4pt}
\begin{picture}(150.00,80.17)(6,0)
\put(80.00,80.00){\makebox(0,0)[cc]{$1$}}
\put(40.00,70.00){\makebox(0,0)[cc]{$1$}}
\put(120.00,70.00){\makebox(0,0)[cc]{$-1$}}
\put(20.00,60.00){\makebox(0,0)[cc]{$\sqrt{2}$}}
\put(60.00,60.00){\makebox(0,0)[cc]{$-\sqrt{2}$}}
\put(120.00,60.00){\makebox(0,0)[cc]{$0$}}
\put(20.00,50.00){\makebox(0,0)[cc]{$0$}}
\put(60.00,50.00){\makebox(0,0)[cc]{$0$}}
\put(100.00,50.00){\makebox(0,0)[cc]{$\sqrt{2}$}}
\put(140.00,50.00){\makebox(0,0)[cc]{$-\sqrt{2}$}}
\put(10.00,40.00){\makebox(0,0)[cc]{$2$}}
\put(30.00,40.00){\makebox(0,0)[cc]{$-2$}}
\put(60.00,40.00){\makebox(0,0)[cc]{$0$}}
\put(100.00,40.00){\makebox(0,0)[cc]{$0$}}
\put(140.00,40.00){\makebox(0,0)[cc]{$0$}}
\put(10.00,30.00){\makebox(0,0)[cc]{$0$}}
\put(30.00,30.00){\makebox(0,0)[cc]{$0$}}
\put(50.00,30.00){\makebox(0,0)[cc]{$2$}}
\put(70.00,30.00){\makebox(0,0)[cc]{$-2$}}
\put(100.00,30.00){\makebox(0,0)[cc]{$0$}}
\put(140.00,30.00){\makebox(0,0)[cc]{$0$}}
\put(10.00,20.00){\makebox(0,0)[cc]{$0$}}
\put(30.00,20.00){\makebox(0,0)[cc]{$0$}}
\put(50.00,20.00){\makebox(0,0)[cc]{$0$}}
\put(70.00,20.00){\makebox(0,0)[cc]{$0$}}
\put(90.00,20.00){\makebox(0,0)[cc]{$2$}}
\put(110.00,20.00){\makebox(0,0)[cc]{$-2$}}
\put(140.00,20.00){\makebox(0,0)[cc]{$0$}}
\put(10.00,10.00){\makebox(0,0)[cc]{$0$}}
\put(30.00,10.00){\makebox(0,0)[cc]{$0$}}
\put(50.00,10.00){\makebox(0,0)[cc]{$0$}}
\put(70.00,10.00){\makebox(0,0)[cc]{$0$}}
\put(90.00,10.00){\makebox(0,0)[cc]{$0$}}
\put(110.00,10.00){\makebox(0,0)[cc]{$0$}}
\put(130.00,10.00){\makebox(0,0)[cc]{$2$}}
\put(150.00,10.00){\makebox(0,0)[cc]{$-2$}}
\put(45.00,75.00){\line(6,1){30.00}}
\put(115.00,75.00){\line(-6,1){31.00}}
\put(23.00,63.00){\line(2,1){13.00}}
\put(57.00,63.00){\line(-2,1){14.00}}
\put(120.00,68.00){\line(0,-1){5.00}}
\put(20.00,58.00){\line(0,-1){5.00}}
\put(60.00,58.00){\line(0,-1){5.00}}
\put(117.00,57.00){\line(-4,-1){14.00}}
\put(123.00,57.00){\line(4,-1){14.00}}
\put(13.00,43.00){\line(4,5){4.00}}
\put(27.00,43.00){\line(-4,5){4.00}}
\put(60.00,48.00){\line(0,-1){5.00}}
\put(100.00,43.00){\line(0,1){5.00}}
\put(140.00,48.00){\line(0,-1){5.00}}
\put(10.00,38.00){\line(0,-1){5.00}}
\put(30.00,33.00){\line(0,1){5.00}}
\put(52.00,34.00){\line(5,4){5.00}}
\put(63.00,38.00){\line(1,-1){4.00}}
\put(100.00,33.00){\line(0,1){5.00}}
\put(140.00,38.00){\line(0,-1){5.00}}
\put(10.00,23.00){\line(0,1){4.00}}
\put(30.00,27.00){\line(0,-1){4.00}}
\put(50.00,23.00){\line(0,1){4.00}}
\put(70.00,27.00){\line(0,-1){4.00}}
\put(92.00,23.00){\line(5,4){5.00}}
\put(104.00,27.00){\line(1,-1){4.00}}
\put(140.00,23.00){\line(0,1){4.00}}
\put(10.00,13.00){\line(0,1){4.00}}
\put(30.00,17.00){\line(0,-1){4.00}}
\put(50.00,13.00){\line(0,1){4.00}}
\put(70.00,17.00){\line(0,-1){4.00}}
\put(90.00,13.00){\line(0,1){4.00}}
\put(110.00,17.00){\line(0,-1){4.00}}
\put(134.00,13.00){\line(4,5){4.00}}
\put(143.00,18.00){\line(4,-5){4.00}}
\put(10.00,7.00){\line(0,-1){2.00}}
\put(30.00,5.00){\line(0,1){2.00}}
\put(50.00,7.00){\line(0,-1){2.00}}
\put(70.00,5.00){\line(0,1){2.00}}
\put(90.00,7.00){\line(0,-1){2.00}}
\put(110.00,5.00){\line(0,1){2.00}}
\put(130.00,7.00){\line(0,-1){2.00}}
\put(150.00,5.00){\line(0,1){2.00}}
\end{picture}
\caption{Kunen's perfect tree \label{kunentree}}
\end{figure}
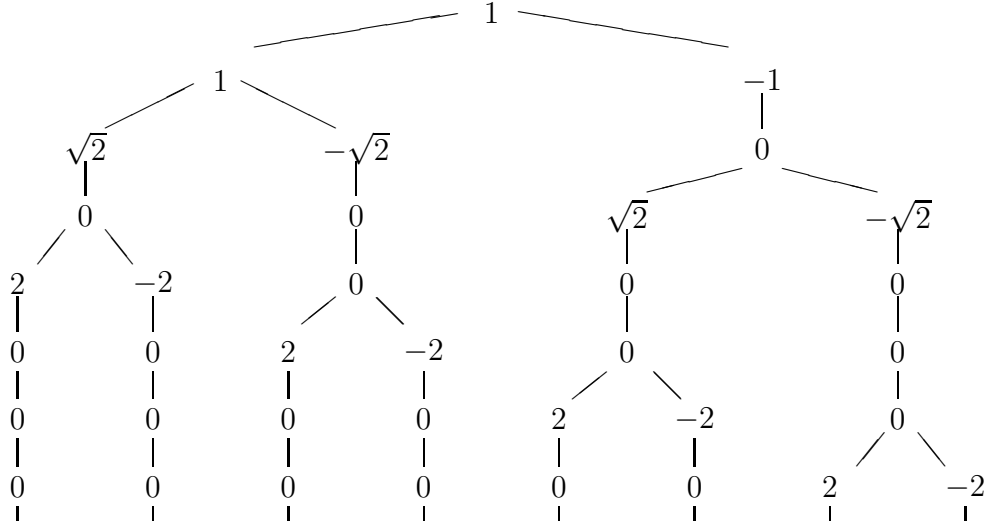

Kunen built a perfect set of pairwise orthogonal elements of
$\rr^\om$ using a different method.  His example is illustrated in
figure \ref{kunentree}.
The next eight levels of his tree use the weights $\sqrt{8}$ and
$-\sqrt{8}$ and so on.

Let $T\su \rr^{<\om}$ be Kunen's tree.  First note that
for any $n\in\om$,
$$T_n=\{s\in T: ||s||=n\}$$
consists of $n$ pairwise orthogonal elements of $\rr^n$.
(This fact was pointed out by Abian.)
Hence each $T_n$ is a maximal orthogonal family in $\rr^n$.
Now suppose $x\in\rr^\om$ is any nontrivial element.
Then there exists $n\in\om$ such that $x\res n$ is nontrivial.
Choose any $s\in T_n$ such that $(s,x\res n )\not=0$.
Suppose for example that $(s,x\res n )>0$ (if its negative a
similar argument works).  Note that for any $t\in T$ either
\begin{itemize}
  \item $t\concat 0\in T$ or
  \item there exists $w>0$ such that
        $t\concat w\in T$ and $t\concat -w\in T$.
\end{itemize}
In other words if a node
doesn't split it continues with 0 and if it does split, then
one way is positive and other negative.  Using this
it is easy to construct a sequence
$$s=s_0\su s_1\su\cdots\su s_m\su\cdots$$
with $s_m\in T_{n+m}$ such that for $m>0$
$$s_m(n+m-1)\cdot x(n+m-1)\geq 0.$$
It follows that
if $b=\cup\{s_m:m\in\om\}$ and $(x,b)$ converges, then
$$(x,b)\geq (x\res n,s)>0,$$
and so, in any case, $x$ is not orthogonal to $b$.
\qed

\pfof{thm2}

For $x\in \rr^\om$ define the support of $x$:
$$\supp(x)=\{n\in\om: x(n)\not=0\}.$$
Let us say that $x$ and $y$ are strongly orthogonal iff
$\supp(x)$ and $\supp(y)$ are almost disjoint
(i.e., intersection is finite) and $(x,y)=0$.
Let $S$ be all the elements of $\rr^\om$ with infinite
support.  Let $S^*$ be all elements $x$ of $S$ such that
if $n\in\supp(x)$  then $|x(n)|>1$.
Let $\po$ be any countable partially ordered set.

\begin{lemma}
Suppose $\{x_n:n\in\om\}\subseteq S^*$ are
pairwise strongly orthogonal. Let
$p\in\po$ and $\tau$ be a $\po$-name such that
$$p\forces \tau\in S \mbox{ and }\forall n \;(x_n,\tau)=0.$$
Then there exists $q\leq p$ and $y\in S^*$ such that
$$q\forces (y,\tau)\not= 0$$
and $y$ is strongly orthogonal to $x_n$ for
every $n\in\om$.
\end{lemma}

\pf

Case 1. $p\forces$``$\supp(\tau)\setminus\cup_{n<N}\supp(x_n)$
is infinite for each $N<\om$''.

\medskip

In this case we will take $q=p$.
Let $\{p_n:n\in\om\}=\{q\in \po:q\leq p\}$ where each element
is listed with infinitely many repetitions.  Build  sequences
$$k_n\in\om\setminus\cup_{m<n}\supp(x_m)$$ and
$$l_n,r_n \in \supp(x_n)\setminus\cup_{m<n}\supp(x_m)$$
such that the sets $\{k_n:n\in\om\}$, $\{l_n:n\in\om\}$,
and $\{r_n:n\in\om\}$ are disjoint and
such that for every $n$
there exists $q\leq p_n$ and $u_n>1$ such that
$$q\forces |u_n\cdot\tau(k_n)|>1.$$
Now construct $y\in S^*$ so that
$$\supp(y)=\{k_n:n\in\om\}\cup\{l_n:n\in\om\}\cup
\{r_n:n\in\om\}$$
and $y(k_n)=u_n$ and $(y,x_n)=0$ for all $n$.
[$(y,x_n)=0$ is accomplished by picking the values of $y(l_n)$ and
$y(r_n)$ inductively.  Use that given $x,u,v\in\rr$ with
$u\not=0$ and $v\not=0$ we can pick $a,b\in\rr$ with $|a|>1$ and $|b|>1$
and $au+bv=x$.  This is possible because if we let $b={x-au \over v}$, then
as $a\to\infty$ we have $|b|\to\infty$.]

But note that
$$p\forces ||\{n: |\tau(k_n)y(k_n)|>1\}||=\om$$
hence we are done with case 1.

Case 2. There exists $r\leq p$, $N<\om$, and $F$ finite such
that
$$r\forces \supp(\tau)\subseteq \cup_{n<N}\supp(x_n)\cup F.$$

In this case find $s\leq r$ and $n_0<N$ such that
$$s\forces ||\supp(\tau)\cap\supp(x_{n_0})||=\om$$
Let $G$ be a $\po$-filter containing $s$.  Let $Q$ be the infinite
set
$$Q=(\supp(x_{n_0})\cap\supp(\tau^G))\setminus
\cup_{n<N, n\not={n_0}}\supp(x_n).$$
Since
$\tau^G$ is orthogonal to $x_{n_0}$ and $|x_{n_0}(n)| >1$ for each
$n\in\supp(x_{n_0})$,
we must be able to find distinct $k_1,k_2\in Q$
such that the vectors $\la\tau^G(k_1),\tau^G(k_2)\ra$ and
$\la x_{n_0}(k_1),x_{n_0}(k_2)\ra$ are not parallel.
[Else
$\la\tau^G(k):k\in Q\ra$ would be parallel to
$\la x_{n_0}(k):k\in Q\ra$ and this would give that
$(x_{n_0},\tau)\not=0$.]
Thus we may find a pair of real numbers $u,v$
(take $\la u,v\ra=\la x_{n_0}(k_2),-x_{n_0}(k_1)\ra$)  such that
$|u|>1$ and $|v|>1$ with $ux_{n_0}(k_1)+vx_{n_0}(k_2)=0$
but $u\tau^G(k_1)+v\tau^G(k_2)\not=0$.  Now build (similarly
to case 1) $y\in S^*$ strongly orthogonal to each $x_n$ and such
that $y(k_1)=u$, $y(k_2)=v$ and
$$\{k_1,k_2\}=\supp(y)\cap (\cup_{n<N}\supp(x_n)\cup F).$$
Let $q\leq s$ be in $G$ such that
$$q\forces u\tau(k_1)+v\tau(k_2)\not=0.$$
Hence $q\forces (y,\tau)\not=0$, so the lemma is proved.

\bigskip

Finally we prove Theorem \ref{thm2}.

In Miller \cite{M}
it is shown that in the Cohen real model the following compactness-like
principle holds:

\begin{quote}
  For any Polish space $X$ and family of Borel sets $\{B_\al: \al<\ka\}$
  where $\omega_1 < \ka < \cc$, if for every $Q\subseteq \ka$  with
  $||Q||<\ka$
  we have $\cap\{B_\al : \al \in Q\}\not=\emptyset$, then
  $\cap\{B_\al : \al <\ka \}\not=\emptyset$.
\end{quote}

It is just stated in \cite{M} for the real line $\rr$ but obviously it
holds for any Borel image of $\rr$.
Note that for any $x \in \rr^\om$, the set
$$B_x=\{y \in \rr^\om\setminus\{{\bf 0}\} : x \mbox { is orthogonal to }y\}$$
is Borel.  Suppose $\{x_\al : \al<\ka\}$ is a pairwise
orthogonal set with $\om_1<\ka<\cc$.  Then for any $Q\subseteq\ka$ of
with $|Q|<\ka$ we have
$$x_\be\in \cap\{B_{x_\al}:\al\in Q\}$$
for any $\be\in\ka\setminus Q$.
By the compactness-like principle
$$\cap\{B_{x_\al}:\al<\ka \}\not=\emptyset$$
and hence $\{x_\al : \al<\ka\}$ is not maximal.

The argument for getting an maximal orthogonal set of
size $\om_1$ is similar to the proof that in the Cohen real
model there is a maximal almost disjoint family of size $\om_1$.
See Kunen \cite{K}.

The problem is to construct a maximal  orthogonal set in
the ground model of CH which is not destroyed by adding
one Cohen real.  To begin with take $\{x_n:n\in\om\}$ in
$S^*$ which are pairwise strongly orthogonal and such that
for every $m<n$ there exists a $k$ with $x_k(m)=1$
and $x_k(i)=0$ for every $i<n$ with $i\not=m$.  These
guarantee that no finite support element of $\rr^\om$ is
orthogonal to every $x_\al$.
Using the continuum hypothesis in the ground model list
pairs $(p_\al,\tau_\al)$ for $\om\leq\al<\om_1$
of elements of $\po \times N$ where $N$ are nice $\po$-names
for potential elements $\rr^\om$.
Build
$x_\al$ for $\al<\om_1$ which are pairwise strongly orthogonal
elements of $S^*$ such that for every $\al\geq\om$ if
$(p_\al,\tau_\al)$  satisfies
that
$$p_\al\forces \tau_\al\in S \mbox{ and }
\forall \be<\al \;(x_\be,\tau_\al)=0,$$
then using the Lemma, there exists $q\leq p_\al$ such that
$$q\forces (x_\al,\tau_\al)\not=0.$$
It follows from this that $\{x_\al:\al<\om_1\}$ remains a
maximal orthogonal set when any number of Cohen reals are added.
This proves Theorem \ref{thm2}.

\pfof{thm3}

The proof uses a modification of a partial order
due to Hechler \cite{H} used to prove the same result for
MAD families.

Let $\rat$ denote the rational numbers.
For $\ga$ an infinite ordinal define $\po_\ga$ as follows.
An element of $\po_\ga$ has the form
$p=(\la s_\al:N\to\rat,P_\al\ra:\al\in F)$
where
\begin{itemize}
  \item $N<\om$ and $F\in [\ga]^{<\om}$,
  \item $s_\al(n)\not=0$ implies $|s_\al(n)|\geq 1$,
  \item (orthogonality) $\al\not=\be\in F$ implies
         $$\sum_{n<N} s_\al(n)s_\be(n)=0,$$
  \item $P_\al\su (F\cap\al)\times N$,
  \item (almost disjoint support) $(\be,n)\in P_\al$ and $n\leq m<N$
  implies   $$s_\al(m)=0 \rmor s_\be(m)=0.$$
\end{itemize}
We define $p\leq q$ iff $F_p\supseteq F_q$, $N^p\geq N^q$, and
$s_\al^p\res{N^q}= s_\al^q$ and $P_\al^p\supseteq P_\al^q$ for every
$\al\in F_q$.  If $G$ is a sufficiently generic $\po_\ga$-filter,
then define
$x_\al:\om\to\rat$ by
$$x_\al(n)=r \rmiff \exists p\in G\; (n<N^p \rmand s^p_\al(n)=r).$$
Note that $(\be,n)\in P_\al$ is a promise that
$\supp(x_\al)\cap\supp(x_\be)\su n$.
Thus the $x_\al$ will be pairwise orthogonal elements of $\rat^\om$ with
almost disjoint support.

\begin{lemma} \label{precond}
  Suppose $p=(\la s_\al:N\to\rat,P_\al\ra:\al\in F)$ is a precondition,
  i.e., it satisfies everything except the orthogonality condition,
  but satisfies instead:
\begin{quote}
   (weak orthogonality) $\al\not=\be\in F$ implies either
         $$\sum_{n<N} s_\al(n)s_\be(n)=0,$$
                        or
   $$\al\notin \dom(P_\be) \rmand \be\notin\dom(P_\al).$$
\end{quote}
  Then $p$ can be extended to a condition in $\po_\ga$.
\end{lemma}
\pf
First list all pairs $\{\la\al_n,\be_n\ra: n<l\}\su [F]^2$ such that
$$\sum_{i<N} s_{\al_n}(i)s_{\be_n}(i)\not=0.$$
Then construct  $s^n_\al:N_n\to\rat$
with
\begin{itemize}
\item $s^p_\al=s^0_\al$ for $\al\in F$,
\item $N_{n+1}=N_n+2$, and
\item $s^{n+1}_\al\res{N_n} = s^n_\al$ for $\al\in F$,
\end{itemize}
as follows. Let
$$x=\sum_{i<N_n}s^n_{\al_n}(i) s^n_{\be_n}(i).$$

Choose $u,v\in \rat$ with $|u|,|v|\geq 1$ and $u+v=-x$.
Now define
$$s^{n+1}_{\al_n}(N_n)=s^{n+1}_{\al_n}(N_{n+1})=1,$$
$$s^{n+1}_{\be_n}(N_n)=u \rmand s^{n+1}_{\be_n}(N_{n+1})=v,$$
and for all other $\de\in F$ define
$$s^{n+1}_{\de}(N_n)=s^{n+1}_{\de}(N_{n+1})=0.$$
Now it is easy to check that
$$(\la s^l_\al:N_l\to\rat,P_\al\ra:\al\in F)$$
is a condition in $\po_\ga$.
\qed

\begin{lemma}
  $\po_\ga$ has ccc, in fact, property K.
\end{lemma}
\pf
Property K means that every uncountable set of conditions contains
an uncountable subset of pairwise compatible conditions.

Given $\Ga$ an uncountable subset of $\po_\ga$
apply a $\De$-system argument
to find $\Si\in[\Ga]^{\om_1}$ and  $F$ such
that $F=F^p\cap F^q=F$
for all distinct $p$ and $q$ in $\Si$.  Next
cutting down $\Si$ we may assume that there exists
$N$ such that $N^p=N$ for all $p\in\Si$ and
there are $(s_\al:\al\in F)$ such that
$$(s^{p}_\al:\al\in F)=(s_\al:\al\in F)$$
for all $p\in \Si$.
Now for any $p,q\in \Si$ define $r$ by:
$$F^r=F^p\cup F^q$$
and
$$\la s_\al^r,P_\al^r\ra =\left\{
\begin{array}{ll}
\la s_\al,P_\al^p\cup P_\al^q\ra   & \mbox{ if $\al\in F$ } \\
\la s_\al^p,P_\al^p\ra             & \mbox{ if $\al\in F^p\se F$ } \\
\la s_\al^q,P_\al^q\ra             & \mbox{ if $\al\in F^q\se F$ } \\
\end{array} \right.$$
The almost disjoint support condition holds for $r$ since it held for
$p$ and $q$.  For
all pairs $\al,\be \in F^r$  $s_\al^r$ and
$s_\be^r$ are orthogonal except possibly those
pairs with $\al\in F^p\se F$ and $\be\in F^q\se F$.  But these
pairs satisfy the weak orthogonality condition and so
by Lemma \ref{precond},  $r$ can be extended to a condition
$\hat{r}$ and clearly $\hat{r}\leq p$ and $\hat{r}\leq q$,
and so $p$ and $q$ are compatible.
\qed

\begin{lemma}
Suppose $G$ is $\po_\ga$-generic over $M$ and $\be<\ga$,
then $G\cap\po_\be$ is $\po_\be$-generic over $M$.
\end{lemma}
\pf

Actually $\po_\be$ is what some author's call a completely
embedded suborder of $\po_\ga$.  This means that for
every $A\su \po_\be$ if $A$ is a maximal antichain of
$\po_\be$, then $A$ is a maximal antichain of
$\po_\ga$.

First note that if $p,q\in \po_\be$ are incompatible
in $\po_\be$, then they are incompatible in $\po_\ga$.
This is because if $r\leq p$ and $r\leq q$ then define
$r\res\be\in \po_\be$ by
$$r\res\be=(\la s_\al^r,P_\al^r\ra:\al\in F^r\cap\be).$$
Then $r\res\be\leq p$ and $r\res\be\leq q$.

Claim. If $p\in \po_\be$ and $r\in \po_\ga$ are incompatible
in $\po_\ga$, then $p$ and  $r\res\be$ are incompatible.

Else suppose there exists $q\in \po_\be$ with $q\leq p$ and $q\leq r\res\be$
and without loss assume $N^q > N^r$.  For
$\al\in F^r\se \be$ define $\hat{s}_\al: N^q\to \rat$ by
$\hat{s}_\al\res{N^r}=s_\al^r$ and $\hat{s}_\al(n)=0$
for all $n$ with $N^r\leq n < N^q$.  Consider the
precondition $t$ defined by
$$F^t=F^r\cup F^q$$
 and
$$\la s_\al^t,P_\al^t\ra =\left\{
\begin{array}{ll}
\la \hat{s}_\al,P_\al^r\ra   & \mbox{ if $\al\in F^r\se\be$ } \\
\la s_\al^q,P_\al^q\ra       & \mbox{ if $\al\in F^q$ }
\end{array} \right.$$
But then $t$ extends to a condition by Lemma \ref{precond},
showing that $p$ and $r$ are compatible and proving the Claim.

It follows from the Claim that if $A\su\po_\be$ is a maximal antichain of
$\po_\be$, then $A$ is a maximal antichain of
$\po_\ga$, and hence the lemma is proved.
\qed

\begin{lemma}  Let $\ga$ be a limit ordinal and \label{finsup}
  suppose $\tau$ is a $\po_\ga$-name for an
  element of $\rr^\om$, $F\in[\ga]^{<\om}$, $H\in\om$,
   and $p\in \po_\ga$
  have the property that
  $$p\forces ||\supp(\tau)||=\om\rmand\supp(\tau)\su
  \bigcup_{\be\in F}\supp(x_\be)\cup H.$$
  Then there exists $\al_1\in\ga$ and $q\leq p$ such that
  $$q\forces (\tau,x_{\al_1})\not=0.$$
\end{lemma}
\pf

Without loss we may assume for every $\al\in F$
$$p\forces (\tau,x_\al)=0.$$
Find $r\leq p$ and $\al_0\in F$ such that
$$r\forces \supp(\tau)\cap \supp(x_{\al_0}) \mbox{ infinite.}$$
Let $G$ be $\po_\ga$-generic with $r\in G$ and
let $Q$ be the infinite set defined by
$$Q=(\supp(\tau^G)\cap \supp(x_{\al_0}))\se (H\cup
\bigcup\{ \supp(x_\be): {\be\in F, \be\not=\al_0}\}).$$
Then there must be  $k_0<k_1\in Q$
such that $\la\tau^G(k_0),\tau^G(k_1)\ra$ is not
parallel to $\la x_{\al_0}(k_0),x_{\al_0}(k_1)\ra$.
Take $t\leq r$ with $t\in G$ with  $H\leq k_0<k_1<N^t-1$
and $s_\al^t(k_0)=u$ and $s_\al^t(k_1)=v$
and such that
$$t\forces v\cdot\tau(k_0) + (-u)\cdot\tau(k_1)\not=0.$$
Note that since $t\in G$, for any $\be\in F\se\{\al_0\}$ we have that
$s^t_\be(k_i)=0$.
Choose any $\al_1> \max(F^t)$ and define the precondition $q$
as follows:
$$F^q=F^t\cup\{\al_1\},$$
$$(\la s^q_\al,P_\al^q\ra : \al\in F^t)=
(\la s^t_\al,P_\al^t\ra : \al\in F^t),$$
and define $s_{\al_1}:N^t\to\rat$
by
$$s_{\al_1}(l)=\left\{
\begin{array}{rl}
v  & \mbox{ if $l=k_0$}\\
-u & \mbox{ if $l=k_1$}\\
0  & \mbox{ otherwise}
\end{array}\right.$$
Define
$$P^q_{\al_1}=\{(0,\be):\be\in F\se \{\al_0\}\}\cup\{(N^t-1,\al_0)\}.$$
The precondition $q$ satisfies all requirements to be an element of
$\po_\ga$ except possibly the orthogonality condition.
But note that for $\be\in F$
$$\sum_{n<N^t} s_{\al_1}(n)s_\be(n)=0.$$
So we only need to worry about $\be\in F^t\se F$ and
$\al_1$.  But for these
$$\al_1\notin \dom(P_\be^q) \rmand \be\notin \dom(P^q_{\al_1})$$
and so we can extend $q$ to a condition using
Lemma \ref{precond}.  Let us denote this extension also by $q$.
By the definition of $P^q_{\al_1}$ we have
that $$q\forces\supp(\tau)\cap\supp(x_{\al_1})=\{k_0,k_1\}$$ and thus
$$q\forces(\tau,x_{\al_1})=v\cdot\tau(k_0) + (-u)\cdot\tau(k_1)\not=0.$$
\qed

\begin{lemma} \label{infsup}
  Suppose $\tau$ is a $\po_\ga$-name for an
  element of $\rr^\om$,
  and let $p\in \po_{\ga+1}$
  have the property that for every $F\in [\ga]^{<\om}$
  $$p\forces ||\supp(\tau)\se
  \bigcup_{\be\in F}\supp(x_\be)||=\om.$$
  Then
  $$p\forces |\tau(n)\cdot x_\ga(n)|\geq 1 \mbox{ for infinitely many $n$}.$$
\end{lemma}
\pf

Suppose not and let $q\leq p$ and $N<\om$ be such
that
$$q\forces \forall n>N\; |\tau(n)\cdot x_\ga(n)|<1.$$
Let $G$ be $\po_{\ga+1}$-generic over $M$ with $q\in G$.
Choose $n>N,N^q$ with
$$n\in\supp(\tau^G)\se\cup \{\supp(x_\al):\al\in F^q\cap\ga\}.$$
Since $\tau$ is a $\po_\ga$-name and $G\cap\po_\ga$ is
$\po_\ga$-generic over $M$ there exists $r\in G\cap\po_\ga$ and
$m\in\om$ such that
$$r\forces |\tau(n)|> {1 \over {m+1}}.$$
We may assume without loss of generality that $N^r>n$ and
$r\leq q\res\ga$.
Note that $s_\al^r(n)=0$ for all $\al\in F^q\cap \ga$.
Now define $\hat{s}_\ga:N^r\to\rat$ as follows.
$$\hat{s}_\ga(k)=\left\{
\begin{array}{ll}
s_\ga^q(k)   & \mbox{ if $k<N^q$ } \\
m+1          & \mbox{ if $k=n$ } \\
0            & \mbox{ otherwise } \\
\end{array} \right.$$
Let $$t=(\la s^r_\al,P^r_\al\ra: \al \in F^r) \cup (\hat{s}_\ga:P^q_\ga).$$
The precondition $t$ satisfies the weak orthogonality condition
of Lemma \ref{precond} and thus can be extended to a condition in
$\po_\ga$.  But it would then force $|x_\ga(n)\cdot \tau(n)|>1$ which
would be a contradiction.
\qed

\begin{lemma}
  Suppose $G$ is $\po_\ga$-generic over $M$ ($\cof(\ga)>\om$)
  and
  $\{x_\al:\al<\ga\}$ are the generic family of mutually
  orthogonal elements of $\rat^\om$, then in $M[G]$ for every
  $y\in \rr^\om$ there exists $\al<\ga$ such that
  $y$ and $x_\al$ are orthogonal.
\end{lemma}
\pf
First note that by any easy density argument for any $m<n<\om$ there
are infinitely many $k<\om$ such that
$$x_k(j)=\left\{
\begin{array}{ll}
1   & \mbox{ if $j=m$ } \\
0   & \mbox{ if $j<n$ and $j\not=m$.}\\
\end{array} \right.$$
Hence we need not worry about $y$ with finite support. If
there exists
$$F\in[\ga]^{<\om} \rmand H<\om$$
such that
$$\supp(y)\su (H\cup\bigcup\{\supp(x_\al):\al\in F\}),$$
then $y$ is taken care of by Lemma \ref{finsup}.  On the other
hand if there is no such $F$, then by using ccc and the fact
that $\cof(\ga)$ is uncountable, we can find $\de<\ga$ and
a $\po_\de$-name $\tau$ for $y$ such that the hypothesis
of Lemma \ref{infsup} holds and thus $(x_\de,y)\not=0$.
\qed

Finally we prove Theorem \ref{thm3}.
Force with the finite support product
$$\sum\{\po_\ga: \ga<\ka\}.$$
This product has property K since each of its factors does.
Since this partial order has cardinality $\ka$ and
$\ka^\om=\ka$ in $M$, the continuum has cardinality $\ka$
in the generic extension.
Also by the product lemma, if $\la G_\ga:\ga<\ka \ra$ is
$\sum\{\po_\ga: \ga<\ka\}$-generic over $M$, then for
each $\ga_0<\ka$ we have that $G_{\ga_0}$ is $\po_{\ga_0}$-generic
over $$M[G_\ga: \ga<\ka, \ga\not=\ga_0].$$  Hence
for each ordinal $\ga<\ka$ of uncountable cofinality
we have a maximal orthogonal family of cardinality $||\ga||$
and Theorem \ref{thm3} is proved.
\qed

\pfof{ma}

For simplicity we first present a proof for the case
when $X$ is disjoint from $l_2$.
Let $\po$ be the following poset.  An element of $\po$ has the
form $$p=(s:N\to\rat, F, P)$$ where $N<\om$, $F\in [X]^{<\om}$, and
$P$ is a finite set of requirements of the form
$P\su F\times (N+1)\times \rat^+$ where $\rat^+$ is the positive rationals and
for every $(x,k,\ep)\in P$,
$$|\sum_{n<k} s(n)\cdot x(n)|<\ep$$
 and for every $l$ with $k<l\leq N$
$$|\sum_{k\leq n <l}s(n)\cdot x(n)|<\ep.$$
We define $p\leq q$ iff $N^q\leq N^p$, $q=p\res{N^q}$, $F^q\su F^p$,
and $P^q\su P^p$.  The poset $\po$ is
$\sigma$-centered, because two conditions with the same $s$ are
compatible.

\begin{lemma} \label{ma1}
  For any $p\in \po$, $x\in F$, and $\ep\in\rat^+$ there exists
  $q\leq p$, and $k<\om$ such that $(x,k,\ep)\in P^q$.
\end{lemma}
\pf

Let $p=(s:N\to\rat,F,P)$.
Let $\ep_0>0$ be such that for any $(y,k,\ga)\in P$
$$|\sum_{k\leq n <N}s(n)\cdot y(n)|+\ep_0<\ga.$$
Choose $N_0\geq N$ so that for any $m>N_0$ and
$y\in F\se \{x\}$
$$|\sum_{N_0\leq n <m}x(n)\cdot y(n)|<\ep_0.$$
Let
$$b= \sum_{n<N} x(n)\cdot s(n).$$
Let $N_1>N_0$ be minimal such that
$$ \sum_{N_0\leq n<N_1} x(n)^2 > b.$$
(This exists since $x$ is not in $l_2$.)   Choose
$\rh$ with $|\rh|\leq 1$ and so that
$$\rh\cdot\sum_{N_0<n<N_1} x(n)^2 = -b.$$
Now consider $t:N_1\to\rr$ defined as below:
$$t(n)=\left\{\begin{array}{cl}
s(n)         & \mbox{ if $n<N$} \\
 0             & \mbox{ if $N\leq n< N_0$} \\
\rh\cdot x(n)  & \mbox{ if $N_0\leq n<N_1$} \\
\end{array}\right.$$

For any $(y,k,\ga)\in P$ with $y\not=x$ and
$m\leq N_1$ note that
$$\begin{array}{lll}
|\sum_{k\leq n <m}t(n) y(n)| & \leq &
|\sum_{k\leq n <N}s(n) y(n)| +
|\rh|\cdot |\sum_{N_0\leq n <m}x(n) y(n)| \\
& \leq &
|\sum_{k\leq n <N}s(n) y(n)| + \ep_0\\
& < & \ga.
\end{array}
$$

Note that since $N_1$ was chosen minimal,
$$|\sum_{k<n<m+1}t(n)\cdot x(n)|\leq |\sum_{k<n<m}t(n)\cdot x(n)|$$
for any $m$ with $N_0\leq m<N_1$, consequently any
requirements involving $x$ are also satisfied.

Since
$$\sum_{n<N_1}t(n)\cdot x(n)=0$$
we can change the values of $t$ on $[N_0,N_1)$ to be rational
to get $t^q:N_1\to\rat$ so that
$$|\sum_{n<N_1}t^q(n)\cdot x(n)|<\ep$$
and still satisfy all the requirements of $P$.
Letting $$q=(t^q,F,P\cup\{(x,N_1,\ep)\})$$ proves the lemma.
\qed

\begin{lemma}\label{ma2}
  For any $p\in\po$ and $l<\om$,
  there exist $q\leq p$
  such that $$\sum_{n<N^q}s^q(n)^2>l.$$
\end{lemma}
\pf
Given $p$ just let $x\in X\se F^p$.  As in the
proof of Lemma \ref{ma1} we can extend $p$ to equal
$x$ as much as we like.  Since $x$ is not in $l_2$
and $p$ has no requirements mentioning $x$ we can
get the norm of $s^q$ greater than $l$.
\qed

The Lemmas show that if we define
$$D_l=\{q\in \po: \sum_{n<N^q}s^q(n)^2>l\}$$
and for each $x\in X$ and $\ep\in\rat^+$
$$D_x^{\ep}=\{q\in \po: x\in F^q \rmand
 \exists k \;(x,k,\ep)\in q\}$$
then these sets are dense.   Applying MA we get
a $\po$-filter $G$ meeting them all.  Then
letting
$$z=\cup\{s^p:p\in G\}$$
we have that $z$ is not in $l_2$ but is orthogonal
to every element of $X$.
This proves Theorem \ref{ma} in the case that $X$ contains
no elements of $l_2$.

\medskip

Next we indicate how to modify our partial order in case
$X$ contains finitely many elements of $l_2$.
Let $H=X\cap l_2$ be finite.
First we replace $\rat$ by any countable field $\rat^*$ which contains
the rationals and all the coefficients of elements of $H$.
We make the following two additional
demands for $$p=(s:N\to\rat^*,F,P)$$ to be an element of
$\po$.
\begin{enumerate}
  \item $H\su F$, and
  \item for each $x\in H$ we have $\sum_{n<N}s(n)x(n)=0$.
\end{enumerate}
The side requirements $P$ are as before including the ones
for elements of $H$.   Now we prove Lemma \ref{ma1} using
our modified definition of $\po$.   In the case that
$x\in H$ this is trivial since we can simply add
$(x,N,\ep)$ to $P^p$.  So let us assume $x\in F\se H$ where
$p=(s:N\to\rat^*,F,P)$.   Let $\ep_0$ be such that
$0<\ep_0<\ep$ and for any $(y,k,\ga)\in P$
$$|\sum_{k\leq n <N}s(n)\cdot y(n)|+\ep_0<\ga.$$
And this time choose $N_0\geq N$ so that for any $m>N_0$ and
$y\in F\se \{x\}$
$$|\sum_{N_0\leq n <m}x(n)\cdot y(n)|<\ep_0/2.$$

For simplicity we begin by giving the proof in the case $H$ has
a single element say
$H=\{z\}$.  If  $\supp(z)$ is finite, it is easy to do
since we can just extend $s$ by zero until we are beyond the
support of $z$ and then apply the same argument as before.
So assume that the support of $z$ is infinite.  Let
$j> N_0$ be so that $z(j)\not=0$.  Choose $\de>0$ so
that
$$\de\cdot\max\{|y(j)|:y\in F\}\leq \ep_0/2.$$
Now choose $N_1>j>N_0$ so that for every $m>N_1$
$$|\sum_{N_1\leq n <m}x(n)z(n)|< |z(j)|\cdot\de.$$
As in the first proof we may find $N_2>N_1$ and $|\rh|\leq 1$
so that
$$\sum_{n<N}x(n)s(n)+\sum_{N_1\leq n<N_2}\rh\cdot x(n) x(n)=0.$$
Now we define $t:N_2\to \rr$ as below:
$$t(n)=\left\{\begin{array}{cl}
 s(n)         & \mbox{ if $n<N$} \\
 0             & \mbox{ if $N\leq n< N_1$ and $n\not=j$} \\
{-1\over z(j)}\sum_{N_1\leq k<N_2}\rh\; x(k)z(k)  & \mbox{ if $n=j$} \\
\rh\cdot x(n)  & \mbox{ if $N_1\leq n<N_2$} \\
\end{array}\right.$$
The value of $t(j)$ is picked to make
$t$ and $z\res{N_2}$ orthogonal (remembering that
$s$ and $z\res N$ are already orthogonal).
First note that
$|y(j)t(j)|<\ep_0/2$ for every $y\in F$, because
$|t(j)|<\de$ and $|y(j)|\cdot\de\leq \ep_0/2$.
Since
$$|\sum_{N_1\leq n < m} y(n)t(n)|<\ep_0/2$$
for all $y\in F\se\{x\}$ and $m\leq N_2$, all the
requirements in $P$ are kept.  Note that
$$|(t,x\res{N_2})|=|t(j)x(j)|\leq \ep_0/2.$$
Now we set $P^t=P\cup \{(x,N_2,\ep)\}$ and $F^t=F$.
The final step is to change the values of $t$ on
$j$ and $[N_1,N_2)$ to elements of $\rat^*$.  First
slightly perturb the values on $[N_1,N_2)$ and then
use them to set the value of $t(j)$ so that
$$t(j)z(j)=-\sum_{N_1\leq n<N_2} t(n)z(n).$$
This new $t$ will be orthogonal to $z\res{N_2}$ and
satisfy all the requirements of $P$.
This concludes the proof of Lemma \ref{ma1} in the
case that $X\cap l_2$ consists of a singleton $\{z\}$.

\medskip

Finally we sketch the proof of the Lemma in the case
that $X\cap l_2=H$ is an arbitrary finite set.
$N_0$ and $\ep_0$ are chosen as before.  Begin
by choosing $H_0\su H$ so that
$\{z\res{[N_0,\om)}: z\in H_0\}$ are linearly independent
(in the space $\rr^{[N_0,\om)}$) and so that
$\{z\res{[N_0,\om)}: z\in H\}$ is contained in the span of
$\{z\res{[N_0,\om)}: z\in H_0\}$.  Choose $N_1>N_0$ so
that $V=\{z\res{[N_0,N_1)}: z\in H_0\}$ are
linearly independent.

\medskip

Claim.  For any $\ep_1>0$ there exists $\de>0$ so
that for any $\la \be_v:v\in V\ra$ with $|\be_v|<\de$
there exists $t\in\rr^{[N_0,N_1)}$ with $||t||_2<\ep_1$
and $(v,t)=\be_v$ for every $v\in V$.

Here  $||\cdot||_2$ is the usual $l_2$ norm
and $(v,t)$ the usual inner product in the finite
dimensional vector space $\rr^{[N_0,N_1)}$.

The proof of the Claim is an elementary exercise in Linear
Algebra.  Let $V=\{v_1,\ldots,v_m\}$ and let $W$ be
the span of $V$.
Define the linear map
$$T:W\to \rr^m
\mbox{ by }
T(t)=\la (t,v_i):i=1,\ldots,m\ra.$$
It follows from the Gram-Schmidt orthogonalization
process that the kernel of $T$ is trivial.  Hence the range
of $T$ is  $\rr^m$.  The existence of $\de$ now follows
from the continuity of $T$.
This proves the Claim.

\medskip

We leave the value of $\ep_1$ to be determined latter.
We find $N_2>N_1$ so that for any $m>N_2$ and
$y\in F\se \{x\}$
$$|\sum_{N_0\leq n <m}x(n)\cdot y(n)|<\de.$$
As in the argument before we find $N_3>N_2$ and
define $t\res{[N_2,N_3)}$ to be a small scaler multiple
of $x\res{[N_2,N_3)}$ in such a way as to make the inner
product of $x$ and $t$ zero.  We now use the
Claim to extend $t$ on the interval $[N_0,N_1)$ to
make sure that for every $y\in H_0$ the inner
product of $t$ and $y$ is zero.  If we choose
$\ep_1$ small enough so that
$$\ep_1\cdot ||x\res{[N_1,N_2)}||_2<\ep_0/2$$
for every $x\in F$,  then this $t$ works (when jiggled
to have range $\rat^*$).

Since $||t\res{[N_1,N_2)}||_2 \cdot ||x\res{[N_1,N_2)}||_2\leq \ep_0/2$
for every $x\in F$ all commitments from
$P$ are honored.  Also $t$ is orthogonal to all elements
of $H$.  The reason is that
$t\res{N_0}=s$ is already orthogonal to every
$z\in H$ and every element of $\{z\res{[N_0,\om)}: z\in H\}$
is contained in the span of
$\{z\res{[N_0,\om)}: z\in H_0\}$ and $t$ is
orthogonal to everything
in $\{z\res{[N_0,\om)}: z\in H_0\}$.

This concludes the proof of Theorem \ref{ma}.

\pfof{support}

For any $h<\om$ define the partial order
$\po_{h}$ as follows.
A condition in $\po_{h}$ has the following form:
$$p=((s_i:N\to\rat^{\not=0}:i<h), R)$$
where $N<\om$, $\rat^{\not=0}$ are the nonzero rationals,
and $R\su [h]^2\times N\times\rat^+$ is a finite set
of requirements satisfying the following:
for any $(\{i,j\},k,\ep)\in R$:
$$|\sum_{n<k} s_i(n)\cdot s_j(n)|<\ep$$
and for any $l$ with $k<l\leq N$
$$|\sum_{k\leq n<l}s_i(n)\cdot s_j(n)|<\ep.$$

\begin{lemma} \label{sup1}
  For any $p\in \po_{h}$, $\ep\in\rat^+$, and
  $\{i,j\}\in [h]^2$ there
  exists $q\leq p$ and $K$ such that
  $(\{i,j\},K,\ep)\in R^q$.
\end{lemma}
\pf
Let
$$b=\sum_{n<N^p}s_i^p(n)\cdot s_j^p(n).$$
Let $N^q=N^p+1$
Define $s_i^q(N^p)=b$ and $s_i^q(N^p)=-1$.  Note
that for any requirement of the form
$((\{i,j\},k,\de)\in R^p$ that
$$\sum_{k\leq n < N^q}s^q_i(n)\cdot s^q_j(n)=
 -b+ \sum_{k\leq n < N^p}s^p_i(n)\cdot s^p_j(n) =
 -\sum_{n < k}s^p_i(n)\cdot s^p_j(n)$$
so the requirement is met.  Now for $l$ different from
$i$ and $j$ we
$s_l(N^p)$ to be some sufficiently small positive rational
which is picked so as to still
satisfy all the requirements of $R^p$.
Let
$$q=((s^q_i:N^p+1\to \rat^{\not=0}:i<h),
R^p\cup\{(\{i,j\},N^p+1,\ep)\}).$$
\qed

\begin{lemma}\label{sup2}
For any $h<\om$ there exists $N<\om$ and
$r_i:N\to\{-1,1\}$ for $i<h$ such that
$r_i$ and $r_j$ are orthogonal for any $i\not=j$.
\end{lemma}
\pf
Let $N=2^{h}$ and identify it with the maps from $h$ into
$2=\{0,1\}$.  For each $i<h$ define $r_i:2^{h}\to \{-1,1\}$ by
$$r_i(t)=\left\{\begin{array}{rl}
     -1 & \mbox{ if $t(i)=1$}\\
     1 & \mbox{ if $t(i)=0$} \\
\end{array}\right.$$
Note that

\centerline{$r_i(t)\cdot r_j(t)=-1$ iff ($t(i)=1$ and $t(j)=0$)
or ($t(i)=0$ and $t(j)=1$).}

But
$$||\{t\in 2^h:t(i)=1 \rmand t(j)=0\}||={1\over 4} N$$
and
$$||\{t\in 2^h:t(i)=0 \rmand t(j)=1\}||={1\over 4} N$$
and so the inner product of $r_i$ and $r_j$ is 0.
\qed

\begin{lemma}\label{sup3}
Suppose $p\in \po_{h}$ and $l<\om$.  Then there
exists $q\leq p$ and $\ep>0$ such that
$(N^q-N^p)\ep^l>1$ and for every
$n$ with $N^p\leq n < N^q$ and $i<h$ we have
$|s^q_i(n)|=\ep$.
\end{lemma}
\pf
Let $\de\in\rat^+$ be such that
for every $(\{i,j\},k,\ga)\in R^p$ we have
$$|\sum_{k\leq n<N^p} s_i^p(n)\cdot s_j^p(n)|+\de^2 <\ga.$$
Apply Lemma \ref{sup2} to obtain pairwise
orthogonal $r_i:N\to\{-1,1\}$ for $i<h$.
Take $\ep={\de \over N}$.  Now define
$\hat{s_i}:N^p+N\to \rat$  by
$$\hat{s_i}(N^p_k)(n)\left\{\begin{array}{ll}
s_i^p(n)        & \mbox{if $n<N^p$}   \\
\ep\cdot r_i(k) & \mbox{if $n=N^p+k$} \\
\end{array}\right.$$
Note that
for any $i\not= j$ and $k<N^p\leq m\leq N^p+N$ that
$$|\sum_{k<n<m}\hat{s}_i(n)\cdot \hat{s}_j(n)|\leq
|\sum_{k<n<N^p}{s}_i^p(n)\cdot {s}^p_j(n)|+\ep^2\cdot N$$
and since $\ep^2\cdot N<\de^2$ any requirements
involving $\{i,j\}$ are kept.
Note also that
$$\sum_{n<N^P+N}\hat{s}_i(n)\cdot \hat{s}_j(n)=
\sum_{n<N}{s}_i(n)\cdot {s}_j(n)$$
and so this trick can be repeated with the same $\ep$
as many times as is necessary to make $(N^q-N^p)\cdot \ep^l>1$.
\qed

Remark. We could have avoided the use of Lemma \ref{sup2}
by taking a pair
of weights $\ep>0$ and $\de>0$ and the column vectors of the
$h \times h$ matrix
with $-\de$ on the diagonal and $\ep$ off the diagonal.
The columns will be orthogonal provided
$$-2\ep\de +(h-2)\ep^2=0 \rmor \de={{h-2}\over 2}\ep$$
In this case we would get that $|s^q_i(n)|\geq\ep$ in
the conclusion of Lemma \ref{sup3}.

\bigskip

Finally, to prove Theorem \ref{support} we construct
a sequence $p_n \in \po_{h_n}$
where $h_n < h_{n+1}$.  Start with any $p_0$ and $h_0$.
At stage $n$ given $p_n$ apply Lemma \ref{sup1} ($h_n$ times)
to obtain $p\leq p_n$ so that for each $\{i,j\}\in [h_n]^2$
for some $k<N^p$ we have that $(\{i,j\},k,{1\over n})\in R^p$.
Now apply Lemma \ref{sup3} to obtain $q\leq p$ and $\ep>0$ such
that $(N^q-N^p)\ep^n>1$ and for every
$m$ with $N^p\leq m < N^q$ and $i<h_n$ we have
$|s^q_i(m)|=\ep$.   Finally obtain $p_{n+1}$ by ``doubling''
$p_n$, i.e., let $h_{n+1}=2\cdot h_n$ define
$$s_{2i}^{p_{n+1}}=s_{2i+1}^{p_{n+1}}=s_i^{p_n}$$
and
$$R^{p_{n+1}}=\{(\{2i,2j\},k,\ga),(\{2i+1,2j+1\},k,\ga):
(\{i,j\},k,\ga)\in R^{p_n}\}.$$
This completes the construction.\footnote{It is not necessary to make
$s_{2i}$ differ from $s_{2i+1}$ since this will
automatically be taken care of when we use of Lemma \ref{sup1}.}
Now define $P$ by
$$P=\{x\in \rr^\om: \forall n \;\exists i<h_n \;x\res{h_n}=s_i^{p_n}\}.$$
It is easy to verify that
the perfect set $P$ has the properties required.
This proves Theorem \ref{support}.
\qed
\pfof{l2}

For part (a) first consider the following example:
$$\begin{array}{rrrrrrrrrrr}
x_0 &=&  &(&1,&-1,& 0, &0,&0,&\ldots&)\\
x_1 &=&  &(&1, &1,&-2, &0,&0,&\ldots&)\\
x_2 &=&  &(&1, &1,& 1,&-3,&0,&\ldots&)\\
\vdots & & & & & & & & &\\
\end{array}
$$
Then $\{x_n: n<\om\}$ is an  orthogonal family in
$l_2$.  We claim that the only $u$'s which are orthogonal
to all of the $x_n$'s are scalar multiples of $(1,1,1,\ldots)$.
If $u=(u_0,u_1,u_2,\ldots)$ is orthogonal to all of the
$x_n$, then
\begin{eqnarray*}
u_0-u_1&=&0\\
u_0+u_1-2u_2&=&0\\
u_0+u_1+u_2-3u_3&=&0\\
\vdots& & \\
\end{eqnarray*}

It follows easily $u_0=u_1=u_2=\cdots$.  This is Kunen's example and
it answered
K.P Hart's original question.

To get the example of part (a)
do as follows.  Work in $\rr^{\om\times\om}$.  Define
$x^n_m$ as follows:
$$x^n_m(l,j)=\left\{\begin{array}{cl}
 1 &  \mbox{if $n=l$ and $j\leq m$}  \\
-(m+1) &  \mbox{if $n=1$ and $j=m+1$}  \\
 0 &  \mbox{otherwise}  \\
\end{array}\right.
$$
So it is the example above
repeated on infinitely many disjoint copies of $\om$.
Let $X=\{x^n_m:n,m<\om\}$.
Then as in Kunen's example:

\centerline{$z\in \rr^{\om\times\om}$ is
orthogonal to all elements of $\{x^n_m:m<\om\}$}
\centerline{iff}
\centerline {$z\res{\{n\}\times\om}$ is a scalar
multiple of the constant 1 sequence.}
Define $y_n:\om\times\om\to\{0,1\}$
to be constantly 1 on $\{n\}\times\om$ and to be 0 otherwise.
Then it is easy to check that $X\cup\{y_n:n<\om\}$ is a maximal
orthogonal family.  Similarly for any $n$ define
$v_n:\om\times\om\to\{0,1\}$
to be constantly 1 on $(\om\se n)\times\om$ and to be 0 otherwise.
Then $X\cup\{y_i:i<n\}\cup\{v_n\}$ is a maximal orthogonal family.
The fact that every orthogonal family extending $X$ is countable
follows easily from noting that if $u,v$ are nontrivial scalar
multiples of $(1,1,1,\ldots)$, then $u$ and $v$ cannot be orthogonal.
This proves part (a).

\bigskip

To prove part (b) we construct a perfect tree $T\su \rr^{<\om}$ as
follows.  Let $T_0=\{()\}$, $T_1=\{(1)\}$ and $T_2=\{(1,1),(1,-1)\}$.
We will do the construction of Kunen (Theorem \ref{thm1-5}) except
we will
use unequal size weights at each stage.  Each $T_n\su\rr^n$ will
consist of $n$ pairwise orthogonal elements of $\rr^n$. The
set $T_{n+1}$
will be obtained by picking exactly one element $s_n$ of $T_n$ and
a pair $\de_n,b_n>0$ with the property that
$(s_n,s_n)=\de_n\cdot b_n$ and then
letting
$$T_{n+1}=\{t\concat 0: t\in T_n, t\not=s_n\}\cup
\{s_n\concat\de_n,s_n\concat -b_n\}.$$
Taking $P=\{x\in\rr^\om:\forall n\; x\res n\in T_n\}$,  the
$s_n$ are chosen so that $P$ has no isolated points and
hence is perfect.  The only remaining things to be picked are the
$\de_n$ and $b_n$. Let $\de_1=b_1=1$.

Given $T_n\su\rr^n$ for $n\geq 2$ let
$$m=\min\{\max\{|(s,x)|:s\in T_n\}: x\in\rr^n,\;
 {3\over 4}\leq ||x||_2\leq 1\}.$$
Here $(s,x)$ refers to the ordinary inner product in $\rr^n$
and $||x||_2$ the $l_2$ norm of $x$.
By compactness it is clear that $m>0$, so we can let
$$\de_n={1\over 2}\cdot\min\{{m\over 2},\de_{n-1}\}$$
and then choose $b_n$ so that $(s_n,s_n)=\de_n\cdot b_n$.

Just as in the proof of Theorem \ref{thm1-5} the family
$P$ is a maximal orthogonal family in $\rr^\om$.
Now define $E\su P$ as follows:
$$E=\{x\in P:\exists n\forall m>n \;(x(m)\in \{0,\de_m\})\}.$$
Since $\de_{m+1}\leq {1\over 2^m}$, it is clear that $E\su l_2$.
We claim $E$ is a maximal orthogonal family in $l_2$.
Suppose $x\in l_2$ is nontrivial and by taking a scalar multiple
(if necessary) assume $||x||_2=1$.  We will find an element of
$E$ which has nonzero inner product with $x$.  Choose $n$
sufficiently large so that $$||x\res n||_2\geq {3\over 4}.$$
Choose $s\in T_n$ so that $\de_n\leq {1\over 4}|(s,x\res n)|$.
Let $y\in E$ be defined by $y\res n=s$ and $y(m)\in\{0,\de_m\}$
for all $m\geq n$.  Note that since $||x||_2=1$, $|x(m)|\leq 1$ and
thus by our choice of $\de_m$'s ($\de_m\leq{1\over 2^{m-n}}\de_m$)
$$|(x\res n, y\res n)| > \sum_{m\geq n} |x(m)y(m)|$$
and consequently $(x,y)\not=0$.
This proves Theorem \ref{l2}.
\qed


\begin{center}
  Open Questions.
\end{center}

\begin{enumerate}

  \item (Abian) Does there exists a model of ZFC with no maximal
  orthogonal family of cardinality $\om_1$?
  In particular, what happens under MA or PFA?

  \item Is it always the case that for any uncountable $\ka$ there
  is a maximal orthogonal family of cardinality $\ka$ iff
  there exists a maximal almost disjoint family of subsets
  of $\om$ of cardinality $\ka$?

  \item (Kunen) If there is a maximal orthogonal family of cardinality
  $\ka$, then does there exists one of cardinality $\ka$ with
  almost disjoint supports?

\end{enumerate}

\bigskip

\begin{center}
  Addresses
\end{center}

\begin{flushleft}
 Arnold W. Miller                             \\
 University of Wisconsin-Madison                 \\
 Department of Mathematics  Van Vleck Hall        \\
 480 Lincoln Drive                                 \\
 Madison, Wisconsin 53706-1388, USA                 \\
 e-mail: miller@math.wisc.edu                        \\

\end{flushleft}

\medskip

\begin{flushleft}
 Juris Steprans                            \\
 York University                            \\
 Department of Mathematics                   \\
 North York,  Ontario M3J 1P3, Canada         \\
 e-mail: juris.steprans@mathstat.yorku.ca     \\
\end{flushleft}

\bigskip
\begin{flushright}
   July 1995.                                \\
\end{flushright}

\end{document}